\documentclass[a4paper, 10pt, reqno]{amsart}
\usepackage{latexsym}
\usepackage{fancyhdr}
\usepackage{amssymb}
\usepackage{array}
\usepackage{enumerate}
\usepackage{verbatim}
\usepackage{xspace}
\usepackage{exscale}
\usepackage[ansinew]{inputenc}
\usepackage{mathrsfs}
\usepackage[all]{xy}
\usepackage{tabularx}
\usepackage[final]{graphicx}

\def\nmb#1#2{#2}         
\newcommand{\subsect}[1]{\par\medskip\noindent{\bf #1.}}
\newenvironment{proclaim}[1]{\par\medskip\noindent{\bf #1.}\it}{\par\smallskip}
\newenvironment{Rem}[1]{\par\medskip\noindent{\bf #1.}}{\par\smallskip}

\numberwithin{equation}{section}

\def\al{\alpha}

\def\ga{\gamma}

\def\ph{\phi}

\def\R{\mathbb{R}}

\let\on=\operatorname
\def\sr#1%
{\ifmmode{}^\dagger\else${}^\dagger$\fi\ifvmode
\vbox to 0pt{\vss
 \hbox to 0pt{\hskip\hsize\hskip1em
 \vbox{\hsize3cm\eightpoint\raggedright\pretolerance10000
 \noindent #1\hfill}\hss}\vss}\else
 \vadjust{\vbox to0pt{\vss%
 \hbox to 0pt{\hskip\hsize\hskip1em%
 \vbox{\hsize3cm\eightpoint\raggedright\pretolerance10000%
 \noindent #1\hfill}\hss}\vss}}\fi%
}

\def\<{\langle}
\def\>{\rangle}
\def\reg{\on{reg}}

\title{Orbit projections as fibrations}
 
\author{Armin Rainer}

\address{Armin Rainer: Fakult\"at f\"ur Mathematik, Universit\"at Wien,
Nordbergstrasse~15, A-1090 Wien, Austria}

\email{armin.rainer@univie.ac.at}

\thanks{}
\date{October 17, 2006}

\begin{document}

\thanks{The author was supported by 
`Fonds zur F\"orderung der wissenschaftlichen Forschung, Projekt P 17108 N04'}
\date{\today}

\maketitle

\vspace{-0.8cm}

\begin{abstract}
The orbit projection $\pi : M \to M/G$ of a proper $G$-manifold $M$ is a fibration if and only 
if all points in $M$ are regular.
Under additional assumptions we show that $\pi$ is a quasifibration if and only if all points are regular.
We get a full answer in the equivariant category: $\pi$ is a $G$-quasifibration if and only if all points are regular.
\end{abstract}

\section{Introduction}

A continuous map is called Hurewicz (Serre) fibration if it has the homotopy lifting property for all topological 
spaces (CW complexes). We show in section \ref{secf} that the orbit projection $\pi : M \to M/G$ of a proper smooth 
$G$-manifold $M$ is a Hurewicz (Serre) fibration if and only if all points in $M$ are regular. 
The proof basically uses the existence of slices at any point and the 
fact that the projection $G \to G/H$, for isotropy subgroups $H$, is a fibration. 
Hence, the result generalizes to proper locally smooth 
$G$-spaces $M$. Moreover, it has its analogon in the category of $G$-spaces and $G$-equivariant maps.  

In section \ref{secqf} we investigate when the orbit projection $\pi$ is a quasifibration, i.e., the canonical inclusion of 
each fiber in the corresponding homotopy fiber is a weak homotopy equivalence. We show that under certain conditions 
$\pi$ is a quasifibration if and only if all points of $M$ are regular. In the equivariant category we get a full 
answer: The orbit projection $\pi$ is a $G$-quasifibration if and only if all points of $M$ are regular.  
Its proof uses some deep theorems from equivariant homotopy theory. 

\section{Orbit projections as fibrations} \label{secf}

\subsect{\nmb.{2.1}. Fibrations} (\cite{dugundji})
Let $E$ and $B$ be topological spaces.
A continuous map $p : E \to B$ is called a {\it Hurewicz fibration} if it has the {\it homotopy lifting property}:
For each topological space $X$, each continuous $f : X \times \{0\} \to E$ and each homotopy 
$\ph : X \times I \to B$ of $p \circ f$, there exists a homotopy $\bar \ph$ of $f$ covering $\ph$.
If $p : E \to B$ has the homotopy lifting property for all CW complexes $X$, it is called a {\it Serre fibration}.
The fibration is {\it regular} if $\bar \ph$ can always be selected to be stationary with $\ph$, i.e.,
for each $x \in X$ such that $\ph(x,t)$ is constant as a function of $t$, the function $\bar \ph(x,t)$ is 
constant as well.  
\[
\xymatrix{
X \times \{0\} \ar[d] \ar[rr]^{f} && E \ar[d]^{p} \\
X \times I \ar[rr]_{\ph} \ar@{-->}[rru]^{\bar \ph} && B
}
\]

A locally trivial fiber bundle is a regular Serre fibration; if the base is paracompact then it is a regular 
Hurewicz fibration. 

For a Serre fibration $p : E \to B$ with fiber $F = p^{-1}(b_0)$ and $p(e_0) = b_0$ the following homotopy sequence is exact:
\begin{equation} \label{hs}
\cdots \to \pi_n(F,e_0) \to \pi_{n}(E,e_0) \to 
\pi_{n}(B,b_0) \to \pi_{n-1}(F,e_0) \to \cdots \to \pi_{0}(B,b_0)
\end{equation}

If $p : E \to B$ is a Hurewicz (Serre) fibration and $B$ is path connected (and all fibers are CW complexes), 
then any two fibers belong to the same homotopy type.

\subsect{\nmb.{2.2}. Compact transformation groups}
Let $G$ be a compact Lie group and let $M$ be a $G$-manifold, i.e., $M$ is a paracompact Hausdorff smooth manifold 
and the action $G \times M \to M, (g,x) \mapsto g.x$ is smooth. 
Endow the orbit space $M/G$ with the quotient topology. 
Then $M/G$ is paracompact and Hausdorff, and the orbit projection $\pi : M \to M/G$ is continuous, open, closed, and proper 
(e.g. \cite{bredon}, \cite{duistermaatkolk}).
In particular, it follows that if $M$ is (path) connected, then $M/G$ is (path) connected. 
The orbits $G.x$ which are exactly the fibers of $\pi$ are compact smooth submanifolds of $M$.
A point $x \in M$ is called stationary if $G.x = \{x\}$.

\begin{proclaim}{Proposition}
Let $G$ be a compact Lie group and let $M$ be an path connected $G$-manifold containing a stationary point. 
The orbit projection $\pi : M \to M/G$ is a Hurewicz (Serre) fibration if and only if every point in $M$ is stationary.
\end{proclaim}

\proof
Suppose that $\pi : M \to M/G$ is a Serre fibration. Since $M/G$ is path connected and all fibers are CW complexes, 
all fibers of $\pi$ belong to the same homotopy type. There is one fiber consisting of one point only, 
namely the stationary point in $M$. It follows that all orbits consist of one point only, since all orbits are 
closed manifolds. Hence each point in $M$ is stationary.
\endproof

\begin{proclaim}{Corollary}
Let $\rho : G \to \on{GL}(V)$ be a representation of a compact Lie group $G$. 
The orbit projection $\pi : V \to V/G$ is a Hurewicz (Serre) fibration if and only if $\rho$ is trivial. \qed
\end{proclaim}

\subsect{\nmb.{2.3}. Proper $G$-manifolds} (\cite{duistermaatkolk})
Let $G$ be a Lie group and let $M$ be a proper $G$-manifold, i.e., the mapping $G \times M \to M \times M, (g,x) \mapsto (g.x,x)$ is proper. Examples are $G$-manifolds where $G$ is compact or properly discontinuous actions on manifolds.
Again $M/G$ is paracompact and Hausdorff, and the orbit projection $\pi$ is continuous, open, closed, and proper.
In a proper $G$-manifold all isotropy subgroups $G_x=\{g \in G : g.x=x\}$ are compact.
We denote by $M_{\reg}$ the set of regular points $x$ in $M$, i.e., points $x$ 
allowing an invariant open neighborhood $U$ such that for all $y \in U$ there exists an equivariant map $f : G.x \to G.y$, 
or equivalently, $G_x \subseteq g G_y g^{-1}$ for some $g \in G$. 
Orbits through regular points are said to be of principal orbit type. 
The orbit types in $M$ are the conjugacy classes $(H)$ of the isotropy subgroups $H \subseteq G$. The inclusion relation on the family of 
isotropy subgroups induces a partial ordering $(H) \le (K)$ on the family of orbit types. 
If $M/G$ is connected, then the minimum orbit type is precisely the principal orbit type.  
If $M$ is a Riemannian $G$-manifold, the regular points in $M$ are exactly those whose slice representation 
$G_x \to \on{O}(T_x (G.x)^{\bot})$ is trivial, where $T_x (G.x)^{\bot}$ denotes the orthogonal complement of 
$T_x (G.x)$ in $T_x M$. 
The set $M_{\reg}$ is open and dense in $M$, and $M_{\reg} \to M_{\reg}/G$ is a locally trivial 
fiber bundle.

\begin{proclaim}{Theorem}
Let $M$ be a proper $G$-manifold. 
The orbit projection $\pi : M \to M/G$ is a Hurewicz (Serre) fibration if and only if $M=M_{\on{reg}}$.
\end{proclaim}

\proof
Suppose that $\pi : M \to M/G$ is a Serre fibration. 
There exists a $G$-invariant Riemannian metric making $M$ a proper Riemannian $G$-manifold. 
By the differentiable slice theorem \cite{palais}, 
for each $x \in M$ there exists a slice $S_x$ such that the $G$-invariant neighborhood $G.S_x$ of $x$ 
is $G$-equivariantly diffeomorphic to the crossed product $G \times_{G_x} S_x$. 
It follows that $G.S_x/G \cong S_x/G_x$ is an open neighborhood of $\pi(x)$ in the orbit space $M/G$.
The slice $S_x$ can be chosen to be the 
diffeomorphic image of an open ball around the origin in the vector subspace $T_x(G.x)^\bot$ of $T_x M$. 
Evidently, the restriction $\pi|_{G.S_x} : G.S_x \to G.S_x/G$ is a Serre fibration as well. 

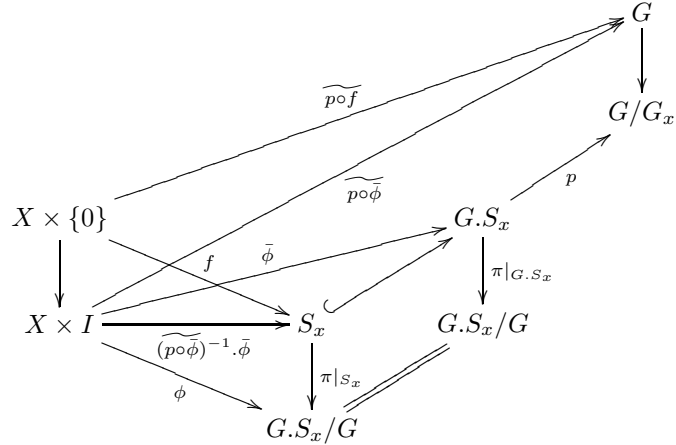
\[
\xymatrix{
&&&& G \ar[d] \\
&&&& G/G_x \\
X \times \{0\} \ar[d] \ar[drr]^(0.57){f} \ar[uurrrr]^{\widetilde{p \circ f}} &&& G.S_x \ar[d]^{\pi|_{G.S_x}} \ar[ur]_{p} & \\
X \times I \ar[drr]_{\ph} \ar[urrr]^{\bar \ph} \ar[uuurrrr]_{\widetilde{p \circ \bar \ph}} 
\ar[rr]_(0.57){(\widetilde{p \circ \bar \ph})^{-1}.\bar \ph} && S_x \ar[d]^{\pi|_{S_x}} \ar@{^(->}[ur] & G.S_x/G & \\
&& G.S_x/G \ar@{=}[ur] &&
}
\]

We claim that also $\pi|_{S_x} : S_x \to G.S_x/G \cong S_x/G_x$ is a Serre fibration. 
Let $f : X \times \{0\} \to S_x$ be continuous and let $\ph : X \times I \to G.S_x/G$ be a homotopy of 
$\pi|_{S_x} \circ f$. Since $\pi|_{G.S_x} : G.S_x \to G.S_x/G$ is a Serre fibration, there exists a 
homotopy $\bar \ph : X \times I \to G.S_x$ of $f$ covering $\ph$. 
Consider the projection $p : G.S_x \cong G \times_{G_x} S_x \to G/G_x$ of the fiber bundle associated to the 
principal bundle $G \to G/G_x$ and the compositions $p \circ f$ and $p \circ \bar \ph$. 
Now $p \circ f$ is constant and equals $eG_x \in G/G_x$ and thus allows a lift into $G$, e.g., $\widetilde{p \circ f}=e$. 
Since $G \to G/G_x$ is a fibration, there exists a homotopy 
$\widetilde{p \circ \bar \ph} : X \times I \to G$ of $\widetilde{p \circ f}$ 
covering $p \circ \bar \ph$. It follows that $(\widetilde{p \circ \bar \ph})^{-1}.\bar \ph : X \times I \to S_x$ 
is a homotopy of $f$ covering $\ph$. Hence the claim is proved.

We may view $\pi|_{S_x} : S_x \to G.S_x/G \cong S_x/G_x$ as the orbit projection of the $G_x$-manifold $S_x$.
Since we may consider the $G_x$-manifold $S_x$ as a 
linear representation of a compact Lie group, the $G_x$-action on $S_x$ must be trivial, by corollary \nmb!{2.2}.
Since $x$ was arbitrary, the statement follows.
\endproof

\begin{Rem}{Remark}
If the orbit projection $\pi : M \to M/G$ is a Hurewicz (Serre) fibration, then it is regular.
\end{Rem}

\subsect{\nmb.{2.4}. $G$-fibrations} (\cite{tomdieck})
A $G$-equivariant continuous map $p : E \to B$ between $G$-spaces $E$ and $B$ is called 
{\it $G$-fibration} if it is a fibration in the category of $G$-spaces and $G$-equivariant maps:
For each $G$-space $X$, each $G$-equivariant continuous $f : X \times \{0\} \to E$ and each $G$-equivariant homotopy 
$\ph : X \times I \to B$ of $p \circ f$, there exists a $G$-equivariant homotopy $\bar \ph$ of $f$ covering $\ph$. 
The $G$-action on $I$ is trivial. It is easy to verify that a $G$-fibration is a fibration in the usual sense. 
More precisely, one can show that $p : E \to B$ is a $G$-fibration if and only if $p^H : E^H \to B^H$ is a 
fibration for each closed subgroup $H \subseteq G$. 
Note that $E^H=\{e \in E : h.e=e ~\text{for all}~h \in H\}$ and $p^H$ denotes the restriction $p|_{E^H}$.

\begin{proclaim}{Theorem}
Let $M$ be a proper $G$-manifold. 
The orbit projection $\pi : M \to M/G$ is a Hurewicz (Serre) $G$-fibration if and only if $M=M_{\on{reg}}$.
\end{proclaim}

\proof
Suppose that $M=M_{\on{reg}}$.
Let $S_x$ be a slice at $x \in M$. Then 
$G.S_x \cong G/G_x \times S_x$, $G.S_x/G \cong S_x$, and
$\pi|_{G.S_x} : G/G_x \times S_x \to S_x$ is given by $([g],s) \mapsto s$. Hence $\pi|_{G.S_x}$ is obviously a
$G$-fibration; a $G$-equivariant homotopy of $f$ covering $\ph$ is given by $\bar \ph = (\on{pr}_{G/G_x} \circ f,\ph)$. 
Since $M/G$ is paracompact and Hausdorff, one can then show that $\pi$ is a $G$-fibration 
analogously with Hurewicz's uniformization theorem (\cite{dugundji}).

Since each $G$-fibration is a fibration, the other implication is an immediate consequence of theorem \nmb!{2.3}.
\endproof

\subsect{\nmb.{2.5}. Proper locally smooth actions}
(\cite{bredon})
Let $G$ be a Lie group and $M$ a proper $G$-space, i.e., $M$ is paracompact and Hausdorff and the $G$-action is continuous. 
Let $G.x$ be an orbit in $M$ and let $V$ be a Euclidean 
vector space on which $G_x$ operates orthogonally. Then a \emph{linear tube} about $G.x$ in $M$ is a 
$G$-equivariant embedding onto an open neighborhood of $G.x$ of the form $G \times_{G_x} V \to M$. 
A $G$-space $M$ is called \emph{locally smooth} if there exists a linear tube about each orbit. 
It that case $M$ must be a topological manifold. It follows from the differentiable slice theorem \cite{palais} 
that proper $G$-manifolds in the sense of \nmb!{2.3} are locally smooth. 

The definition of local smoothness can by extended to manifolds $M$ with boundary. For this we require, for 
orbits $G.x$ lying in the boundary of $M$, tubes of the form $G \times_{G_x} V^+ \to M$, where 
$V^+ = \{y \in \R^n : y_1 \ge 0\}$ and $G_x$ acts orthogonally on $V^+$ (in particular the $y_1$-axis is stationary).  

Properness guarantees that all isotropy groups $G_x$ are compact. 
Hence the orbits in each $G_x$-space $V$ (resp. $V^+$) are compact 
manifolds and therefore 
CW complexes. 
It follows that the arguments in \nmb!{2.2}, \nmb!{2.3}, and \nmb!{2.4} are applicable and we obtain 

\begin{proclaim}{Theorem}
Let $M$ be a proper locally smooth $G$-space (with boundary). 
The orbit projection $\pi : M \to M/G$ is a Hurewicz (Serre) ($G$-)fibration if and only if $M=M_{\on{reg}}$. \qed
\end{proclaim}

\section{Orbit projections as quasifibrations} \label{secqf}

\subsect{\nmb.{3.1}. Quasifibrations}
(\cite{hatcher})
A continuous map $p : E \to B$ with $B$ path connected is called \emph{quasifibration} if the induced map 
$p_* : \pi_n(E,p^{-1}(b),e) \to \pi_n(B,b)$ is an isomorphism for all $b \in B$, $e \in p^{-1}(b)$, and $n \ge 0$, 
or equivalently, if the inclusion of each fiber $p^{-1}(b)$ into the homotopy fiber $F_b$ of $p$ over $b$ 
is a weak homotopy equivalence. The fiber $p^{-1}(b)$ is included in 
\[
F_b = \{(e,\ga) \in E \times C^0(I,B) : \ga(0)=p(e), \ga(1)=b\}
\] 
as the pairs $(e,\ga)$ with $e \in p^{-1}(b)$ and 
$\ga$ the constant path at $b$.  
If $B$ is not path connected, then $p : E \to B$ is a \emph{quasifibration} if the restriction of $p$ over each 
path component of $B$ is a quasifibration. 
For quasifibrations the homotopy sequence \eqref{hs} is exact. 
All fibers of a quasifibration $p : E \to B$ with $B$ path connected belong to the same weak homotopy type.
Hurewicz and Serre fibrations are quasifibrations.

\begin{proclaim}{\nmb.{3.2}. Lemma}
Let $M$ be a proper $G$-manifold with connected orbit space $M/G$. Let $k$ be the least number of connected components of 
isotropy groups of dimension $m := \on{min} \{\on{dim} G_x : x \in M\}$. Then the following conditions are equivalent:
\begin{enumerate}
\item[$(1)$] $G.x$ is a principal orbit.
\item[$(2)$] $G_x$ is of dimension $m$ and has $k$ connected components. 
\end{enumerate}
\end{proclaim}

\proof
Let $x \in M$ such that $G_x$ has minimal dimension and the least number of connected 
components for this dimension in all of $M$. Let $S_x$ be a slice at $x$. For any $y \in G.S_x$ we have 
$y \in g.S_x=S_{g.x}$ and thus $G_y \subseteq G_{g.x} = g G_x g^{-1}$, for some $g \in G$. 
By the choice of $x$, we find $G_y = g G_x g^{-1}$ which shows that $G.x$ is principal.
The converse implication follows from the fact that there is precisely one principal orbit type, if $M/G$ is connected. 
\endproof

\begin{proclaim}{\nmb.{3.3}. Theorem}
Let $M$ be a proper $G$-manifold.
Let one of the following conditions be satisfied:
\begin{enumerate}
\item[$(1)$] $G$ is finite.
\item[$(2)$] $G$ is compact, connected, and simply connected.
\item[$(3)$] $G$ is compact and there exists a connected and simply connected orbit in each path component of $M/G$.
\item[$(4)$] There exists a weakly contractible orbit in each path component of $M/G$.
\end{enumerate}
Then the orbit projection $\pi : M \to M/G$ is a quasifibration if and only if $M=M_{\on{reg}}$.
\end{proclaim}

\proof
Suppose that $\pi : M \to M/G$ is a quasifibration.
We may suppose that $M/G$ is path connected, by restricting $\pi$ over each path component of $M/G$
and treating them separately. Then all orbits belong to the same weak homotopy type. 

We claim that each of the four conditions in the theorem implies that all occurring isotropy groups have the same dimension 
and the same number of connected components.
If $G$ is finite, all orbits and thus all isotropy groups have the same cardinality. 
Assume that $G$ is compact. Let $G.x$ and $G.y$ be distinct orbits. Since they are compact manifolds, 
we find $\on{dim} G.x = \on{dim} G.y$, and, consequently, $\on{dim} G_x = \on{dim} G_y$. If $(2)$ is satisfied we may conclude from 
the homotopy sequences of the fibrations $G \to G/G_x \cong G.x$ and $G \to G/G_y \cong G.y$ that 
$\pi_0(G_x) \cong \pi_1(G/G_x) \cong \pi_1(G/G_y) \cong \pi_0(G_y)$. 
Assume that $(3)$ holds true. Then each orbit is connected and simply connected. Let $G.x$ be principal and $G.y$ arbitrary. 
Without loss $G_x \subseteq G_y$ and we have a locally trivial fiber bundle $G/G_x \to G/G_y$ with fiber $G_y/G_x$. The associated 
homotopy sequence yields that $G_y/G_x$ is trivial, whence the statement.
Finally, if condition $(4)$ is fulfilled, all orbits are weakly contractible, whence any two isotropy groups have the same 
weak homotopy type, again by the homotopy sequence of $G \to G/G_x \cong G.x$. Since all isotropy groups are compact 
manifolds, the claim follows.
  
Since $M/G$ is connected, there is precisely one principal orbit type, namely the type corresponding to the isotropy 
group with minimal dimension and minimal number of connected components. 
By the claim, all points are regular. 
\endproof

\subsect{\nmb.{3.4}. $G$-quasifibrations}
A $G$-equivariant continuous map $p : E \to B$ between $G$-spaces $E$ and $B$ is called 
{\it $G$-quasifibration} if $p^H : E^H \to B^H$ is a quasifibration for each closed subgroup $H \subseteq G$.
In particular, a $G$-quasifibration is a quasifibration.
Any $G$-fibration is a $G$-quasifibration.

\begin{proclaim}{Corollary}
Let $M$ be a proper $G$-manifold. Suppose that one of the conditions $(1)$ -- $(4)$ in the theorem \nmb!{3.3} is satisfied.
Then the orbit projection $\pi : M \to M/G$ is a $G$-quasifibration if and only if $M=M_{\on{reg}}$.
\end{proclaim}

\proof
The statement follows from theorem \nmb!{2.4} and theorem \nmb!{3.3}.
\endproof

\subsect{\nmb.{3.5}}
Let $M$ be a proper $G$-manifold. For the orbit projection $\pi : M \to M/G$ consider the path fibration 
$p : E_\pi \to M/G$, where 
\[
E_\pi = \{(x,\ga) \in M \times C^0(I,M/G) : \ga(0) = \pi(x)\}
\] 
and $p(x,\ga) = \ga(1)$. 
The space $C^0(I,M/G)$ carries the compact-open topology and $E_\pi$ inherits the subspace topology from $M \times C^0(I,M/G)$.
Then $p : E_\pi \to M/G$ is a Hurewicz fibration with fibers 
\[
F_z= \{(x,\ga) \in M \times C^0(I,M/G) : \ga(0) = \pi(x), \ga(1) = z\}, 
\]
$z \in M/G$.
The $G$-action on $M$ induces a natural $G$-action on each fiber $\pi^{-1}(z)$ and each homotopy fiber $F_z$ for which 
the canonical inclusion $\pi^{-1}(z) \hookrightarrow F_z$ is equivariant.

Assume that $M/G$ is path connected.
Let $\pi^{-1}(u) = G.x$ and $\pi^{-1}(z) = G.y$ be distinct orbits in $M$ and choose a path $\al : I \to M/G$ with 
$\al(0)=u$ and $\al(1)=z$. We have the following diagram 
\[
\xymatrix{
\pi^{-1}(u) \ar@{^(->}[r] & F_u \ar[r] & F_z & \pi^{-1}(z) \ar@{_(->}[l]
},
\]
where $F_u \to F_z$ given by $(x,\ga) \mapsto (x,\al \ga)$ is a homotopy equivalence. Note that each arrow in the diagram is 
$G$-equivariant.

\begin{proclaim}{Theorem}
Let $M$ be a proper $G$-manifold and assume that $M/G$ is path connected. The following conditions are equivalent:
\begin{enumerate}
\item[$(1)$] The orbit projection $\pi : M \to M/G$ is a $G$-quasifibration.
\item[$(2)$] The inclusion $\pi^{-1}(z) \hookrightarrow F_z$ is a homotopy equivalence allowing a $G$-equivariant 
homotopy inverse for all $z \in M/G$.
\item[$(3)$] Let $(H)$ be the principal orbit type. For any orbit type $(K)$ there is a (weak) homotopy 
equivalence $f : G/H \to G/K$, and $f$ is $G$-equivariant. 
\item[$(4)$] $M = M_{\on{\reg}}$.
\end{enumerate}
\end{proclaim}

\proof
It is evident that $(2)$ implies $(3)$.

Let us assume that condition $(3)$ is satisfied. We prove $(4)$. Without loss we may suppose that $H \subseteq K$ and have the commuting diagram
\[
\xymatrix{
H \ar@{^(->}[r] \ar@{^(->}[d] & G \ar[r] \ar[d]^{\on{id}} & G/H \ar[d]^{f} \\ 
K \ar@{^(->}[r] & G \ar[r] & G/K
}
\] 
Consequently, using the fact that $G \to G/H$ and $G \to G/K$ are fibrations, we obtain the commuting diagram
\[
\xymatrix{ 
\pi_{n+1}(G) \ar[r] \ar[d] & \pi_{n+1}(G/H) \ar[r] \ar[d] & \pi_{n}(H) \ar[r] \ar[d] & \pi_{n}(G) \ar[r] \ar[d] & \pi_{n}(G/H) \ar[d] \\ 
\pi_{n+1}(G) \ar[r] & \pi_{n+1}(G/K) \ar[r] & \pi_{n}(K) \ar[r] & \pi_{n}(G) \ar[r] & \pi_{n}(G/K)
}
\]
for each $n \ge 0$ where the rows are exact and all vertical arrows apart from the middle one are isomorphisms. 
By the five lemma, the vertical middle arrow is an isomorphism as well. It follows that $\pi_n(H) \cong \pi_n(K)$ for all $n \ge 0$, 
and, by Whitehead's theorem, we find that $H$ and $K$ are homotopically equivalent. Since $H$ and $K$ are compact Lie groups, 
we may conclude that they have the same dimension and the same number of connected components. By lemma \nmb!{3.2}, 
all points in $M$ have to be regular.  

Theorem \nmb!{2.4} yields that $(4)$ implies $(1)$.


Finally, we prove that $(1)$ implies $(2)$. If $\pi : M \to M/G$ is a $G$-quasifibration, then the fixed point maps 
$\pi^H : M^H \to M/G$ are quasifibrations for all closed subgroups $H \subseteq G$. Let $z \in M/G$ be arbitrary. 
The fiber $\pi^{-1}(z)$ and the homotopy fiber $F_z$ are $G$-spaces in a canonical way, and we have 
$(\pi^H)^{-1}(z) = \pi^{-1}(z)^H$ and 
\[
F_z^H = \{(x,\ga) \in M^H \times C^0(I,M/G) : \ga(0)=\pi^H(x),\ga(1)=z\}.
\]  
Since $\pi^H$ is a quasifibration, the canonical inclusion $\pi^{-1}(z)^H \hookrightarrow F_z^H$ is a weak 
homotopy equivalence. By \nmb!{3.6.1}, we may conclude that the $F_z^H$ are homotopy 
equivalent to CW complexes and that $F_z$ has the $G$-homotopy type of a $G$-CW complex. 
By Whitehead's theorem, the inclusions $\pi^{-1}(z)^H \hookrightarrow F_z^H$ are homotopy equivalences. 
By \nmb!{3.6.2}, we obtain that the inclusion $\pi^{-1}(z) \hookrightarrow F_z$ is even a $G$-homotopy equivalence. 
Hence $(2)$. 
\endproof

\begin{proclaim}{Corollary}
Let $M$ be a proper locally smooth connected $G$-space (with boundary) and suppose that the principal orbits are of 
codimension $1$. 
Then the orbit projection $\pi : M \to M/G$ is a quasifibration if and only if $M=M_{\on{reg}}$.
\end{proclaim}

\proof
It is proved in \cite{bredon} IV.8. that under these condition either all orbits are principal or $M$ is equivalent as 
$G$-space to the mapping cylinder of the equivariant map $G/H \to G/K$ for $(H)$ principal and $(H) < (K)$ or to the union of 
the two mapping cylinders of $G/H \to G/K_i$ for $(H)$ principal and $(H) < (K_i)$, $i=0,1$. 
In the latter cases the orbit space $M/G$ is isomorphic either to $[0,1)$ or to $[0,1]$, and the natural projection of 
a mapping cylinder identifies with $\pi$. This projection is a quasifibration if and only if the mapping inducing the 
mapping cylinder is a weak homotopy equivalence. By the implication $(3) \Rightarrow (4)$ in the forgoing theorem, 
the statement of the corollary follows.
\endproof

\subsect{\nmb.{3.6}. Equivariant homotopy theory}
We collect a few results from equivariant homotopy theory needed in the proof of theorem \nmb!{3.5}.

For a definition of $G$-CW complexes see the cited references.

\begin{proclaim}{\nmb.{3.6.1}. Result}
(\cite{waner} 4.14; \cite{maysigurdsson} 3.3.5; see also \cite{milnor})
Suppose that $X$ and $Y$ are proper $G$-spaces.
Let $f : X \to Y$ be a $G$-map, and let $y \in Y$ have isotropy group $H$. Then, regarding $f$ as an $H$-equivariant 
map based at $y$, the homotopy fiber $F_y$ of $f$ has the $H$-homotopy type of an $H$-CW complex whenever $X$ and $Y$ 
have the $G$-homotopy type of $G$-CW complexes. 
\end{proclaim}

It is proved in \cite{illman} that a proper $G$-manifold $M$ has a $G$-CW structure.
The orbit space $M/G$ is triangulable by \cite{verona}.

\begin{proclaim}{\nmb.{3.6.2}. Result}
(\cite{lueck} 1.1, \cite{tomdieck} II.2.7)
A $G$-map $f : X \to Y$ of $G$-CW complexes is a $G$-homotopy equivalence if and only if for any subgroup
$H \subseteq G$ which occurs as isotropy subgroup of $X$ or $Y$ the induced map $f^H : X^H \to Y^H$ is a 
homotopy equivalence.
\end{proclaim}

\end{document}